\documentclass[a4paper,10pt]{article}

\usepackage{amsfonts,amssymb,amsmath,xypic}

\numberwithin{equation}{section}

\usepackage{enumitem}

\usepackage[normalem]{ulem}

\usepackage{tikz-cd}
\usepackage{comment}
\newtheorem{theo}[equation]{Theorem}
\newtheorem{coro}[equation]{Corollary}
\newtheorem{lemm}[equation]{Lemma}
\newtheorem{prop}[equation]{Proposition}
\newtheorem{defi}[equation]{Definition}
\newtheorem{rema}[equation]{Remark}

\newtheorem{exam}[equation]{Example}

\newenvironment{proof}{\noindent \textbf{{Proof.}} \sf}

\def\qed{\hfill $\square$ \bigskip}

\def\s{{\mathcal S}}

\def\lim{\mathop{\rm lim}\nolimits}

\def\DD{\mathsf D}
\def\EE{\mathsf E}

\def\HH{\mathsf H}
\def\HHH{\mathsf{HH}}

\def\Ext{\mathsf{Ext}}
\def\Hom{\mathsf{Hom}}
\def\Tor{\mathsf{Tor}}

\def\Ker{\mathsf{Ker}}

\def\Coker{\mathsf{Coker}}

\def\Im{\mathsf{Im}}

\def\dim{\mathsf{dim}}

\def\tauL{\HHH^*_\tau(\Lambda)}
\def\tauLn{\HHH^n_\tau(\Lambda)}
\def\tauXn{\HH^n_\tau(\Lambda, X)}
\def\tauXi{\HH^i_\tau(\Lambda, X)}

\def\sfTr{\mathsf{Tr}}
\def\tauhhX{\HH_n^\tau(\Lambda,X)}
\def\tauhhL{\HHH_n^\tau(\Lambda)}

\def\tauhhLs{\HHH_*^\tau(\Lambda)}

\begin{document}
\sf

\title{Happel's question, Han's conjecture and $\tau$-Hochschild (co)homology}
\author{Claude Cibils,  Marcelo Lanzilotta, Eduardo N. Marcos,\\and Andrea Solotar
\thanks{\tiny This work has been supported by the projects PIP-CONICET 11220200101855CO, USP-COFECUB.
The third mentioned author was supported by the thematic project of FAPESP  2018/23690-6,  research grants from CNPq 308706/2021-8  and  310651/2022-0. The fourth mentioned author is a research member of CONICET (Argentina), Senior Associate at ICTP and visiting Professor at Guangdong Technion-Israel Institute of Technology.}}

\date{}
\maketitle
\begin{abstract}
We introduce the $\tau$-Hochschild (co)homology of a finite dimensional associative algebra  $\Lambda$ by means of the higher Auslander-Reiten translate of O. Iyama.  We show that the global dimension  of $\Lambda$,  Happel's question and Han's conjecture are   related to the $\tau$-Hochschild (co)homology.
\end{abstract}



%\tableofcontents

\section{\sf Introduction}

Let $k$ be a field and $\Lambda$ a finite dimensional associative $k$-algebra with Jacobson radical $r$, such that  $E=\Lambda/r$ is separable. For short,  such algebras  will be just called \emph{algebras}.  All modules and bimodules we consider are finite dimensional. We denote $\otimes_k$ by $\otimes$ and we use the symbol = whenever a canonical isomorphism exists. 

The contents of this paper are as follows. In Section \ref{tau HH (co)homology} we define $\tau$-Hochschild (co)homology in positive degrees using the higher Auslander-Reiten translate considered by O. Iyama in \cite{IYAMA2007, IYAMA2007b}. The definition stems from one of the main ideas of $\tau$-tilting theory, see \cite{ADACHIIYAMAREITEN2014, TREFFINGER 2023, MARSH2023}. 
 We also prove that $\tau$-Hochschild (co)homology is Morita invariant. However it is not derived invariant, see Example \ref{tau not derived}.

In general the dimensions of the $\tau$-Hochschild (co)homology are strictly greater in each degree than the corresponding ones for Hochschild (co)homology. This is shown by the computations of $\tau$-Hochschild (co)homology of radical square algebras. We postpone these calculations to the last section \ref{rad^2=0} to ensure continuity in the development of the theory. 

 Section \ref{Happel's} contains known results that we will use later. First we recall  Happel's result on the minimal projective resolution of an algebra in \cite{HAPPEL1989}. We also record that if the global dimension of an algebra is $d$, then the Hochschild (co)homology is zero in degrees  greater than $d$.  Y. Han and B. Keller in \cite{HAN2006, KELLER1998} proved that for these algebras, actually the Hochschild homology is zero in positive degrees.

In Section \ref{dimensions} we first prove that if an algebra is of finite global dimension $d$, then the $\tau$-Hochschild (co)homology is zero in degrees greater than $d$. Moreover the $\tau$-Hochschild homology is also zero in degree $d$, but in general it is not zero in degree $d-1$, as Example \ref{no go further} shows. However we prove that for bound quiver algebras whose quivers have no oriented cycles, the $\tau$-Hochschild homology is $0$.

We also provide formulas for the dimension of the $\tau$-Hochschild (co)ho\-mo\-lo\-gy in degree $n$ with coefficients in a bimodule, see Theorem  \ref{dim tau HH}. The formulas involve the dimensions of the  Hochschild (co)ho\-mo\-lo\-gy in degrees strictly smaller than $n$, and the dimensions of the Tor (or Ext) vector spaces of simple modules multiplied by the dimension of the corresponding isotypic components of the bimodule. 

We call a graded vector space $V_*$ \emph{infinite} if for infinitely many $n$ we have $V_n\neq 0$. Otherwise, we call it \emph{finite}. 

In Section \ref{tau Happel's and Han's},  we define the following: a bound quiver algebra $\Lambda$ has \emph{infinite + global dimension} (resp. of \emph{infinite co+ global dimension}) if there exists a pair of vertices $(x,y)$ of the quiver such that 
\begin{itemize}
   \item $y\Lambda x\neq 0$,
   \item $\Tor^\Lambda_*(k_x,{}_yk)$ (resp. $\Tor^\Lambda_*(k_y,{}_xk)$) is infinite
 \end{itemize}
where $k_y$ is the simple right $\Lambda$-module associated to the vertex $y$, and ${}_xk$ is the simple left $\Lambda$-module associated to the vertex $x$.
Clearly, if $\Lambda$ is of infinite + and/or  infinite  co+ global dimension then $\Lambda$ is of infinite global dimension. We do not know counterexamples for the converse statement.

 D. Happel  proved in  \cite{HAPPEL1989} that if the Hochschild cohomology of an algebra is infinite, then its global dimension is infinite. He wrote in \cite[p. 110]{HAPPEL1989}  ``The converse seems to be not known". This phrase would latter become known as ``Happel's question". Commutative algebras are positive answers, see \cite{AVRAMOVIYENGAR 2005}, but the family of local algebras considered in \cite{BUCHWEITZGREENMADSENSOLBERG2005} are negative answers to Happel's question: they have Hochschild cohomology zero in degrees greater than or equal to $3$, however they are of infinite global dimension - as all non trivial local algebras. In contrast their $\tau$-Hochschild cohomology is infinite, see Example \ref{BGMS}. 
 
 A main result of this paper is that a bound quiver algebra is of infinite co+ global dimension if and only if its $\tau$-Hochschild cohomology is infinite.
 
 Y. Han conjectured in \cite{HAN2006} that if the global dimension of an algebra is infinite, then its Hochschild homology is infinite. This has been proved for several families of algebras, see for instance \cite{AVRAMOVVIGUE,  BERGHERDMANN, BERGHMADSEN2009, BERGHMADSEN2017, BACH, HAN2006, SOLOTARSUAREZVIVAS, SOLOTARVIGUE}.
 
Another main result of this paper is that a bound quiver algebra is of infinite + global dimension if and only if its $\tau$-Hochschild homology is infinite.
 
In Section \ref{examples of + and co+} we examine the possibility that for a bound quiver algebra, infinite global dimension could imply infinite + and  infinite  co+ global dimension. 

First, any non trivial local algebra is of infinite + and  infinite  co+ global dimension. Second, let $\EE (\Lambda)$ be the Yoneda $k$-category of  an algebra $\Lambda$. Note that $\EE(\Lambda)$ is infinite dimensional if and only if $\Lambda$ is of infinite global dimension. We prove that if $\EE(\Lambda)$ admits a $k$-subcategory which is infinite dimensional although finitely generated, then $\Lambda$ is of infinite + and  infinite co+ global dimension. Third, algebras with many non zero Peirce components - for the precise statement see Proposition \ref{todos no cero} - verify the above possible implication.  Fourth, bound quiver algebras with a loop in their quiver and verifying the extension conjecture, also satisfy the above possible implication. Moreover in Subsection \ref{?}, we record that if a bound quiver algebra of infinite global dimension were not of infinite + global dimension, then it will disprove Han’s conjecture. If a bound quiver algebra of infinite global dimension were not of infinite co+ global dimension, then it would  be a negative answer to Happel's question.

\textbf{Acknowledgements}:  We thank Hip\'olito Treffinger for comments regarding our definition of $\tau$-Hochschild (co)homology and the higher Auslander-Reiten translate considered by O. Iyama.

\section{\sf $\tau$-Hochschild (co)homology}\label{tau HH (co)homology}

Let $A$ be an algebra and $M$, $N$ left $A$-modules.  Let $\DD=\Hom_k(-,k)$ and $\sfTr$ denote the transpose, see for instance \cite{AUSLANDERREITENSMALO1995}. Recall that the Auslander-Reiten translate is $\tau=\DD\sfTr$, see \cite{AUSLANDERREITEN1975,AUSLANDERREITENSMALO1995}.

Following a main idea of $\tau$-tilting theory in \cite{ADACHIIYAMAREITEN2014, TREFFINGER 2023, MARSH2023},  whenever $\Ext^1_A(M,N)$ appears we replace it with $\DD\Hom_A(N,\tau M)$, as we have done in \cite{CIBILSLANZILOTTAMARCOSSOLOTAR2024}. For the main definitions and properties  of $\tau$-tilting theory, see for instance \cite{IYAMAREITEN2014} or \cite{TREFFINGER 2023}.

This idea is based on the Auslander-Reiten duality formula, see \cite{AUSLANDERREITEN1975}:
$$\Ext^1_A(M,N) = \DD\overline{\Hom}_A(N,\tau M),$$
where $\overline{\Hom}_A(N,\tau M)$ is the quotient of  $\Hom_A(N,\tau M)$ by the subspace of morphisms which factor through injective modules. Replacing $\overline{\Hom}_A(N,\tau M)$ by \\ $\Hom_A(N,\tau M)$ can be interpreted as recovering those morphisms.

Let $\Lambda$ be an algebra and let $X$ be a $\Lambda$-bimodule, considered as a left module over the enveloping algebra $\Lambda \otimes \Lambda^{\mathsf{op}}$.  The Hochschild cohomology and homology in degree $n\geq  0$  are respectively (see \cite{HOCHSCHILD1945}):
$$\HH^n(\Lambda, X) = \Ext^n_{\Lambda \otimes \Lambda^{\mathsf{op}}} (\Lambda, X),$$
$$\HH_n(\Lambda, X) = \Tor_n^{\Lambda \otimes \Lambda^{\mathsf{op}}} ( X, \Lambda).$$
 As usual, we will denote $\HHH^n(\Lambda) =\HH^n(\Lambda,\Lambda)$ and $\HHH_n(\Lambda) =\HH_n(\Lambda,\Lambda)$. 

From now on we will replace $\Lambda \otimes \Lambda^{\mathsf{op}}$ by $\Lambda - \Lambda$, since left $\Lambda \otimes \Lambda^{\mathsf{op}}$-modules are the same as $\Lambda$-bimodules.

Next consider Heller's syzygy functors $\{\Omega^n\}$ for $\Lambda$-bimodules, see \cite{HELLER}. In what follows  $\tau$ stands for the Auslander-Reiten translate for $\Lambda$-bimodules.  
The following definition is due to O. Iyama in \cite[p. 56]{IYAMA2007b} and \cite{IYAMA2007}. See also for instance \cite{JASSOKULSHAMMER2079,JASSOKVAMME2019}.

\begin{defi}
Let $n\geq 1$. The \emph{$n$-Auslander-Reiten translate} is
$$\tau_n = \tau \Omega^{n-1}.$$
\end{defi}

\begin{prop}\label{classicalHH}
For all $n\ge 1$
$$\HH^n(\Lambda, X)= \DD\overline{\Hom}_{\Lambda - \Lambda}(X,\tau_n\Lambda).$$
$$\HH_n(\Lambda, X)= \overline{\Hom}_{\Lambda - \Lambda}(\DD X,\tau_n\Lambda).$$
\end{prop}

\begin{proof}
\begin{align*}
\bullet\  \HH^n(\Lambda, X) &=\Ext^n_{\Lambda - \Lambda} (\Lambda, X)= \Ext^1_{\Lambda - \Lambda}(\Omega^{n-1}\Lambda, X) =  \DD\overline{\Hom}_{\Lambda - \Lambda}(X,\tau \Omega^{n-1}\Lambda).\\
\bullet\   \HH_n(\Lambda, X) &=\Tor_n^{\Lambda - \Lambda} (X,\Lambda)=\DD\Ext^n_{\Lambda - \Lambda} (\Lambda, \DD X) 
 = \DD  \Ext^1_{\Lambda - \Lambda}(\Omega^{n-1}\Lambda, \DD X)\\&=\DD\DD \overline{\Hom}_{\Lambda - \Lambda}(\DD X,\tau \Omega^{n-1}\Lambda) = \overline{\Hom}_{\Lambda - \Lambda}(\DD X,\tau \Omega^{n-1}\Lambda).
\end{align*}
\qed
\end{proof}
 In view of Proposition \ref{classicalHH} we set the following. 
\begin{defi}
The \emph{$\tau$-Hochschild cohomology and homology} of $\Lambda$ with coefficients in $X$ in degree $n\geq 1$ are respectively

$$\HH^n_\tau(\Lambda,X) = \DD\Hom_{\Lambda - \Lambda}(X, \tau_n \Lambda)$$
$$\tauhhX = \Hom_{\Lambda - \Lambda}(\DD X, \tau_n \Lambda).$$
 We denote  $\HHH^n_\tau(\Lambda) =\HH^n_\tau(\Lambda,\Lambda)$ and $\HHH_n^\tau(\Lambda) =\HH_n^\tau(\Lambda,\Lambda).$
\end{defi}

As for Hochschild (co)homology, we have the following

\begin{prop}\label{double duality}For all $n\ge 1$
$$\DD\ \tauhhX = \HH^n_\tau(\Lambda,\DD X).$$
\end{prop}
\begin{proof}
$$\DD\HH^n_\tau(\Lambda,\DD X) = \DD \DD\Hom_{\Lambda - \Lambda}(\DD X, \tau_n \Lambda) = \Hom_{\Lambda - \Lambda}(\DD X, \tau_n \Lambda) = \tauhhX.$$\qed
\end{proof}

\begin{rema}

     In \cite{CIBILSLANZILOTTAMARCOSSOLOTAR2024} we introduced $\tau$-Hochschild cohomology  in degree one. The above definition agrees since  $\tau_1=\tau$. 
\end{rema}

\begin{lemm}\label{inclusion}
For all $n\geq 1$, there is an inclusion $$\HH^n(\Lambda, X) \hookrightarrow \tauXn$$ and a surjection $$\tauhhX\twoheadrightarrow \HH_n(\Lambda,X).$$
\end{lemm}
\begin{proof}
The surjection $${\Hom}_{\Lambda - \Lambda}(X, \tau_n \Lambda)\twoheadrightarrow  \overline{\Hom}_{\Lambda - \Lambda}(X, \tau_n \Lambda)$$
gives
  $$\HH^n(\Lambda, X) = \DD\overline{\Hom}_{\Lambda - \Lambda}(X, \tau_n \Lambda)\hookrightarrow \DD{\Hom}_{\Lambda - \Lambda}(X, \tau_n \Lambda)=\tauXn.
  $$

For Hochschild homology the surjection is
$$\tauhhX=\Hom_{\Lambda - \Lambda}(\DD X, \tau_n \Lambda)\twoheadrightarrow\overline{\Hom}_{\Lambda - \Lambda}(\DD X, \tau_n  \Lambda)=\HH_n(\Lambda,X).$$
\qed
\end{proof}

\begin{rema}
In general the dimensions of the $\tau$-Hochschild cohomology (resp. homology) spaces of an algebra are strictly greater than the dimensions of the respective Hochschild cohomology (resp. homology) spaces, see Section \ref{rad^2=0}.
 \end{rema}

 Let $A$ be an algebra and $M$, $N$ be left $A$-modules. The next lemma  is a consequence of the following canonical isomorphisms of vector spaces
$$\DD\Hom_A(N,M)=\DD M\otimes_A N $$
$$\Hom_A(M,N)=\Hom_A(\DD N, \DD M).$$

\begin{lemm}\label{taucoho as transpose of Heller}
For $n\geq 1$ we have
  $$\tauXn= \DD \tau_n \Lambda\otimes_{\Lambda-\Lambda} X,$$
 $$\tauhhX= \Hom_{\Lambda-\Lambda}( \DD \tau_n \Lambda, X).$$
\end{lemm}
\begin{proof}
We have 
\begin{itemize}
  \item $\HH^n_\tau(\Lambda,X) = \DD\Hom_{\Lambda - \Lambda}(X, \tau_n \Lambda)= \DD \tau_n \Lambda\otimes_{\Lambda-\Lambda} X$
  \item $\tauhhX = \Hom_{\Lambda - \Lambda}(\DD X, \tau_n \Lambda)=\Hom_{\Lambda-\Lambda}( \DD \tau_n \Lambda, \DD\DD X).$
\end{itemize}\qed
\end{proof}

\normalsize
Consider the minimal projective resolution of the $\Lambda$-bimodule $\Lambda$
\begin{equation}\label{minimal projective resolution}
\cdots\longrightarrow P_n\stackrel{d_n}{\longrightarrow}P_{n-1}\longrightarrow \cdots \longrightarrow P_2\stackrel{d_2}{\longrightarrow}P_{1}\stackrel{d_1}{\longrightarrow}P_{0}\stackrel{d_0}{\longrightarrow}\Lambda\longrightarrow 0.
\end{equation}

The Hochschild cohomology $\HH^*(\Lambda, X)$ is the cohomology of the following cochain complex:
\begin{align}\label{cochains}
\begin{split}
&0 \longrightarrow \Hom_{\Lambda-\Lambda}(P_0, X)   \stackrel{\delta_1}{\longrightarrow}    \Hom_{\Lambda-\Lambda}(P_1, X)         \stackrel{\delta_2}{\longrightarrow}   \Hom_{\Lambda-\Lambda}(P_2, X)\longrightarrow\\&\cdots \longrightarrow   \Hom_{\Lambda-\Lambda}(P_{n-1}, X)   \stackrel{\delta_n}{\longrightarrow}   \Hom_{\Lambda-\Lambda}(P_n, X) \longrightarrow\cdots
\end{split}
\end{align}
while the Hochschild homology $\HH_*(\Lambda, X)$ is the homology of the following chain complex:
\begin{align}\label{chains}
\begin{split}
& \cdots \longrightarrow X\otimes_{\Lambda-\Lambda}P_n \stackrel{\delta'_n}{\longrightarrow} X\otimes_{\Lambda-\Lambda}P_{n-1} \longrightarrow \cdots \\& \longrightarrow X\otimes_{\Lambda-\Lambda}P_2 \stackrel{\delta'_2}{\longrightarrow}  X\otimes_{\Lambda-\Lambda}P_1 \stackrel{\delta'_1}{\longrightarrow} X\otimes_{\Lambda-\Lambda}P_0 \longrightarrow 0. 
\end{split}
\end{align}

Note that $\HH^*(\Lambda,X)$ and $\HH_*(\Lambda,X)$ might also be computed using any projective resolution.

\begin{theo}\label{coker}
 For $n\geq 1$ we have $$\tauXn= \Coker\ \delta_n ,$$
 $$\tauhhX =\Ker\ \delta'_n.$$

\end{theo}
\begin{proof}
 By definition $\Omega^{n-1}\Lambda =\Ker \ d_{n-2}=\Im \ d_{n-1}$ in (\ref{minimal projective resolution}). Moreover the minimal projective presentation of $\Omega^{n-1}\Lambda$ is:
$$ P_n\stackrel{d_n}{\longrightarrow}P_{n-1}\longrightarrow \Omega^{n-1}\Lambda \longrightarrow 0$$
By definition of the transpose of $\Lambda$-bimodules, the cokernel of  $$d_n^*=\Hom_{\Lambda -\Lambda}(d_n,    \Lambda\otimes\Lambda)$$ is $\sfTr\Omega^{n-1}\Lambda,$ thus it fits into an exact sequence
\begin{equation}\label{Transpose of Omega}
\Hom_{\Lambda -\Lambda}(P_{n-1} ,  \Lambda\otimes\Lambda)  \stackrel{d_n^*}{\longrightarrow}   \Hom_{\Lambda -\Lambda}(P_{n} ,  \Lambda\otimes\Lambda) \longrightarrow  \sfTr\Omega^{n-1}\Lambda \longrightarrow 0.
\end{equation}
The functor $ - \otimes_{\Lambda - \Lambda} X$ is right exact, so we obtain an exact sequence
\begin{align}\label{transpose of Heller}
\begin{split}
  \Hom_{\Lambda -\Lambda}(P_{n-1} ,  \Lambda\otimes\Lambda) \otimes_{\Lambda - \Lambda} X \longrightarrow \  &\Hom_{\Lambda -\Lambda}(P_{n} ,  \Lambda\otimes\Lambda) \otimes_{\Lambda - \Lambda} X \\&\longrightarrow  (\sfTr\Omega^{n-1}\Lambda) \otimes_{\Lambda - \Lambda} X\longrightarrow 0.
  \end{split}
  \end{align}
By Lemma \ref{taucoho as transpose of Heller} we have:
\begin{align}\label{tau coho}
\begin{split}
  \Hom_{\Lambda -\Lambda}(P_{n-1} ,  \Lambda\otimes\Lambda) \otimes_{\Lambda - \Lambda} X \longrightarrow \  &\Hom_{\Lambda -\Lambda}(P_{n} ,  \Lambda\otimes\Lambda) \otimes_{\Lambda - \Lambda} X \\&\longrightarrow  \tauXn  \longrightarrow 0.
  \end{split}
\end{align}
  On the other hand for any algebra $A$, any projective left $A$-module $P$ and any  left $A$-module $M$, the following holds
  \begin{equation}\label{igualdad buscada}
\Hom_A(P,A)\otimes_A M= \Hom_A (P, M).
  \end{equation}
Indeed, this can be verified for $P=A$, then for free modules, and finally for direct summands of free modules.

By using (\ref{igualdad buscada}) for $\Lambda$-bimodules,  (\ref{tau coho}) is isomorphic to the exact sequence
\begin{align*}
 \Hom_{\Lambda-\Lambda}(P_{n-1}, X)   \stackrel{\delta_n}{\longrightarrow}   \Hom_{\Lambda-\Lambda}(P_n, X) \longrightarrow \Coker \ \delta_n \longrightarrow 0.
\end{align*}
Hence $\tauXn= \Coker \ \delta_n.$

The proof of the isomorphism $\tauhhX =\Ker\ \delta'_n$  is analogous; we apply the functor $\Hom_{\Lambda-\Lambda}(-,X)$ to the exact sequence (\ref{Transpose of Omega}), obtaining
\begin{align}\label{tau ho}
\begin{split}
0 \to \tauhhX \to \Hom_{\Lambda -\Lambda}( &\Hom_{\Lambda -\Lambda}(P_{n} ,  \Lambda\otimes\Lambda), X) \to \\ &\Hom_{\Lambda -\Lambda}( \Hom_{\Lambda -\Lambda}(P_{n-1} ,  \Lambda\otimes\Lambda), X) 
  \end{split}
\end{align}

Let $A$ be an algebra, $P$ a projective left $A$-module and $M$ a right $A$-module. As in (\ref{igualdad buscada}) we have
$$\Hom_A(\Hom_A(P,A), M) = M\otimes_A P.$$
Thus we obtain that (\ref{tau ho}) is isomorphic to the exact sequence
\begin{align*}
0 \to \Ker \ \delta'_n \to X\otimes_{\Lambda-\Lambda}P_{n} \stackrel{\delta'_n}{\to} X\otimes_{\Lambda-\Lambda}P_{n-1} 
\end{align*}
\qed
\end{proof}
\begin{rema}
As in Theorem \ref{coker}, consider the minimal projective resolution of an algebra as a bimodule, and the corresponding  complexes of cochains and chains with respect to a bimodule.

Theorem \ref{coker} says that   the difference between the Hochschild cohomology and the $\tau$-Hoch\-schild cohomology is that for the former we compute cocycles modulo coboundaries, while for the latter we compute all the cochains modulo coboundaries.

Analogously, for the $\tau$-Hochschild homology we consider the cycles, but without making the quotient by the boundaries - the latter gives the Hochschild homology.

\end{rema}

The Hochschild (co)homology is derived invariant, however the $\tau$-Hochschild (co)homology is only Morita invariant  as it is shown in the next result. Example \ref{tau not derived} gives two derived equivalent algebras with non isomorphic $\tau$-Hochschild (co)homologies.

\begin{coro}
The $\tau$-Hochschild (co)homology is Morita invariant.
\end{coro}
\begin{proof}
If $\Lambda$ and $\Lambda'$ are Morita equivalent algebras, then their enveloping algebras are Morita equivalent. Moreover, their respective minimal projective resolutions correspond through the equivalence between $\Lambda$ and $\Lambda'$-bimodules. \qed
\end{proof}

\section{\sf Happel's minimal resolution and the $\Tor$ functor}\label{Happel's}
Let $\Lambda$ be an algebra. As stated in the Introduction the Jacobson radical is denoted $r$ and $E=\Lambda/r$ is separable. Note that $E$ is a semisimple algebra and also a semisimple $\Lambda$ and $E$--bimodule. By the Wedderburn-Mal'tsev decomposition (see \cite{WEDDERBURN 1908, MAL'TSEV 1942}), there exists a subalgebra of $\Lambda$ still denoted $E$, such that $\Lambda =E \oplus r$.

Next we explicitly provide the bimodules of the minimal projective resolution (\ref{minimal projective resolution}) in terms of $\Tor_n^\Lambda(E,E)$, see \cite{BUTLER KING 1997}. D. Happel showed in \cite{HAPPEL1989} that the multiplicities of the projective bimodules in this resolution are given in terms of the dimensions of the Ext vector spaces of simple left $\Lambda$-modules.

\begin{rema} \label{ext and tor}
If  $A$, $B$ and $C$ are algebras and $M$  and $N$ are respectively  $B-A$ and $A-C$-bimodules, then $\Tor_n^A(M,N)$ and $\Ext^n_A(\DD M, N)$ are respectively $B-C$ and $C-B$ bimodules. There is a canonical isomorphism of $C-B$-bimodules
$$\DD\Tor_n^A (M,N) = \Ext^n_A(N, \DD M) \mbox{ for } n\geq 0.$$

In particular the $E$-bimodule $\Tor_n^\Lambda (E,E)$ is isomorphic to $\DD \Ext^n_\Lambda(E, \DD E)$. Note that $\DD E$ is isomorphic to $E$ as a $\Lambda$-bimodule. Then $\Tor_n^\Lambda (E,E) = \DD \Ext^n_\Lambda(E, E)$ and 
$$\dim_k \Tor_n^\Lambda (E,E)= \dim_k \Ext^n_\Lambda(E, E).$$
\end{rema}

For the proof of the next result, we first recall the following. Let $\Lambda$ be an algebra with Jacobson radical $r$, a projective left $\Lambda$-module $Q$ and a left $\Lambda$-module $M$.  A surjective morphism $f:Q\to M$ is a \emph{projective cover} if and only if $\Ker f \subset rQ$. Therefore a projective resolution of $M$
$$\cdots\longrightarrow Q_n\stackrel{\Delta_n}{\longrightarrow} Q_{n-1} \longrightarrow \cdots Q_1\stackrel{\Delta_1}{\longrightarrow} Q_0 \stackrel{\Delta_0}{\longrightarrow} M \longrightarrow 0$$
is \emph{minimal} if and only if $\Im \Delta_{n+1}\subset rQ_n$ for $n\geq 0$. Indeed, $\Ker \Delta_n = \Im \Delta_{n+1}$.

\begin{theo}\label{minimal projective resolution Happel}\cite{BUTLER KING 1997,HAPPEL1989}\label{Happel}
 The projective $\Lambda$-bimodule $P_n$ in the minimal projective resolution (\ref{minimal projective resolution})  of $\Lambda$ is isomorphic to
 $$\Lambda\otimes_E\Tor^\Lambda_n(E,E)\otimes_E\Lambda.$$
\end{theo}
\begin{proof}
Any projective $\Lambda$-bimodule is isomorphic to
$\Lambda\otimes_E T\otimes_E\Lambda$ for some $E$-bimodule $T$, so in (\ref{minimal projective resolution}) we write  $P_n= \Lambda\otimes_E T_n\otimes_E\Lambda$, for an $E$-bimodule $T_n$ and for each $n$.

Any projective $\Lambda$-bimodule is also left and right projective. Hence the resolution has a contracting homotopy of right (or left) modules. Let $M$ be a left $\Lambda$-module. Applying the functor $-\otimes_\Lambda M$ to (\ref{minimal projective resolution}) we obtain a projective resolution of $M$
\begin{align}\label{M}
\begin{split}
\cdots &\longrightarrow\Lambda\otimes_E T_{n+1} \otimes_E M \stackrel{\Delta_{n+1}}{\longrightarrow}
\Lambda\otimes_E T_{n}\otimes_E M
\longrightarrow
\cdots\\
&\longrightarrow
\Lambda\otimes_E T_2 \otimes_E M
\stackrel{\Delta_2}{\longrightarrow}
\Lambda\otimes_E T_1 \otimes_E M
\stackrel{\Delta_1}{\longrightarrow}
\Lambda\otimes_E T_0 \otimes_E M
\stackrel{\Delta_0}{\longrightarrow}
M\longrightarrow 0
\end{split}
\end{align}
which in general is not minimal. Our purpose is to prove that for $M=E$ this resolution is minimal. The Jacobson radical of the enveloping algebra $\Lambda \otimes \Lambda^{\mathsf{op}}$ is $r\otimes \Lambda + \Lambda \otimes r$. The minimality of the resolution of $\Lambda$ is equivalent to
$$\Im \ d_{n+1} \subset  (r\otimes_E T_{n}\otimes_E \Lambda) + (\Lambda\otimes_E T_{n} \otimes_E r)$$
for all $n\geq 1.$ 

For $M=E$ we have $rE=0$. Hence the projective resolution (\ref{M}) for $M=E$ is minimal:
\begin{align}\label{E}
\begin{split}
\cdots &\longrightarrow\Lambda\otimes_E T_{n+1} \stackrel{\Delta_{n+1}}{\longrightarrow}
\Lambda\otimes_E T_{n}
\longrightarrow
\cdots\\
&\longrightarrow
\Lambda\otimes_E T_2
\stackrel{\Delta_2}{\longrightarrow}
\Lambda\otimes_E T_1
\stackrel{\Delta_1}{\longrightarrow}
\Lambda\otimes_E T_0
\stackrel{\Delta_0}{\longrightarrow}
E\longrightarrow 0.
\end{split}
\end{align}
Furthermore, applying the functor $E\otimes_\Lambda -$ yields a complex whose homology is $\Tor_n^\Lambda (E,E)$:
$$\cdots \to T_{n+1} \to T_{n} \to \cdots \to T_2 \to T_1 \to T_0 \to 0.$$
Since $Er=0$ the morphisms of the above chain complex are $0$, hence
$$T_n = \Tor^\Lambda_n(E,E).$$\qed

\end{proof}

The following results have been proved by D. Happel in \cite{HAPPEL1989} as immediate consequences of Theorem \ref{minimal projective resolution Happel}.

\begin{coro}\cite{HAPPEL1989}
The global dimension of an algebra $\Lambda$ equals the projective dimension of $\Lambda$ as a $\Lambda$-bimodule.
\end{coro}

\begin{coro}\cite{HAPPEL1989}\label{fgd implies hh finite}
If the algebra $\Lambda$ is of finite global dimension $d$, then for any $\Lambda$-bimodule $X$ and for $n \geq d+1$, we have $\HH^n(\Lambda, X)=0=\HH_n(\Lambda, X)$. 
\end{coro}

Let $Q$ be a \emph{quiver}, that is a finite oriented graph,  with finite set of \emph{vertices} $Q_0$, finite set of \emph{arrows} $Q_1$ and $s,t: Q_1\to Q_0$ the maps giving the \emph{source and target} of each arrow. A \emph{path} $\gamma = \gamma_n\dots \gamma_1$ is a sequence of $n$ concatenated arrows, that is $t(\gamma_i)=s(\gamma_{i+1})$ for $i=1, \dots, n-1$. We set $s(\gamma)=s(\gamma_1)$ and $t(\gamma)=t(\gamma_n)$. 

The \emph{path algebra} $kQ$ is the tensor algebra over $kQ_0$ of the $kQ_0$-bimodule $kQ_1$.

Let $F$ be the ideal spanned by the arrows of $Q$. The quotient algebra $\Lambda = kQ/I$ where the ideal $I$ is admissible is called a \emph{bound quiver algebra}. We have $r= F/I$ and $E=\Lambda/r= kQ_0$.

\begin{rema} \sf Let $\Lambda=kQ/I$ be a bound quiver algebra.
\begin{itemize}
\item For $x\in Q_0$, we have ${}_xk = \DD k_x.$

\item For $x$ and $y$ vertices,  the simple $\Lambda$-bimodule ${}_yk_x$ is $ {}_yk\otimes k_x.$

\item The bimodule $E$ decomposes as $E=\oplus_{x\in Q_0}\ {}_x k_x$.

\item Let $X$ be a $\Lambda$-bimodule. Its $E-E$-isotypic component of type ${}_yk_x$ is $yXx$.
\end{itemize}

\end{rema}

\begin{prop}\label{pedazos de Tor}For all $n\ge 0$
  $$y\Tor_n^\Lambda (E,E)x= \Tor_n^\Lambda (k_y,{}_xk)=\DD \Ext^n_\Lambda({}_xk, {}_yk).$$
\end{prop}
\begin{proof}
Consider the projective resolution (\ref{M}) of ${}_xk$:
\begin{align*}\label{yK}
\begin{split}
\cdots &\longrightarrow\Lambda\otimes_E T_nx \stackrel{d_n}{\longrightarrow}
\Lambda\otimes_E T_{n-1}x
\longrightarrow
\cdots\\
&\longrightarrow
\Lambda\otimes_E T_2x
\stackrel{d_2}{\longrightarrow}
\Lambda\otimes_E T_1x
\stackrel{d_1}{\longrightarrow}
\Lambda\otimes_E T_0x
\stackrel{d_0}{\longrightarrow}
{}_xk\longrightarrow 0.
\end{split}
\end{align*}
After applying $k_y\otimes_\Lambda -$ we obtain a chain complex whose homology is $\Tor_n^\Lambda (k_y,{}_x k)$:
$$\cdots \to yT_nx \to yT_{n-1}x \to \cdots \to yT_2x \to yT_1x \to yT_0x \to 0.$$
The morphisms of this chain complex are $0$ for the same reason as in the proof of Theorem \ref{Happel}. Hence
$$yT_nx = \Tor^\Lambda_n(k_y,{}_xk).$$
By Theorem \ref{Happel} we know that $T_n =  \Tor_n^\Lambda (E,E).$ The last equality of the statement is a consequence of Remark (\ref{ext and tor}).
\qed
\end{proof}

\section{\sf Dimensions of the $\tau$-Hochschild (co)homology}\label{dimensions}

The dimensions of $\tau$-Hochschild cohomology and $\tau$-Hochschild homology are in general strictly greater than the dimensions of  Hochschild cohomology and homology respectively, as shown for instance in Section \ref{rad^2=0}. Despite of that, Corollary \ref{fgd implies hh finite} has an analog as follows.

\begin{prop}\label{tauHH also}
If the algebra $\Lambda$ is of finite global dimension $d$, then for any $\Lambda$-bimodule $X$ and for $n\geq d+1$, we have $\HH^n_\tau(\Lambda, X)=0=\HH_n^\tau(\Lambda, X)$. Also $\HHH_d^\tau(\Lambda)= 0$.
\end{prop}
\begin{proof}
Theorem \ref{minimal projective resolution Happel} ensures that $P_n=0$ for $n\geq d+1$. Thus for $n\geq d+1$   
$$\tauXn= \Coker\left(\Hom_{\Lambda-\Lambda}(P_{n-1}, X)   \stackrel{\delta_n}{\longrightarrow}   \Hom_{\Lambda-\Lambda}(P_n, X)\right) = 0$$
and
 $$\tauhhX =\Ker\left( X\otimes_{\Lambda-\Lambda}P_n       \stackrel{\delta'_n}{\longrightarrow}  X\otimes_{\Lambda-\Lambda}P_{n-1}  \right) =0,$$
see Theorem \ref{coker}.

It remains to prove that $\HHH_d^\tau(\Lambda)=0$.  The chain complex whose homology is $\HHH_d^\tau(\Lambda)$ is 
$$ 0\stackrel{\delta'_{d+1}}{\longrightarrow} \Lambda\otimes_{\Lambda-\Lambda} P_d \stackrel{\delta'_d}{\longrightarrow}  \Lambda\otimes_{\Lambda-\Lambda} P_{d-1} \to \cdots \to  \Lambda\otimes_{\Lambda-\Lambda} P_0 \to 0.$$
Y. Han and B. Keller proved in \cite[Proposition 6]{HAN2006} and \cite{KELLER1998}, that for algebras of finite global dimension the equality $\HHH_n(\Lambda) =0$ holds for $n>0$. Hence 
$$0=\HHH_d(\Lambda) = \frac{\Ker \delta'_d}{\Im \delta'_{d+1}}=\Ker \delta'_d = \HHH_d^\tau(\Lambda).$$
\qed
\end{proof}

\begin{rema}
The Hochschild homology of an algebra of finite global dimension vanishes in positive degrees. This is no longer the case for $\tau$-Hochschild homology, see Example \ref{no go further}. 
However if $\Lambda$ is a bounded quiver algebra whose quiver has no oriented cycles, then $\HHH_n^\tau(\Lambda)= 0$ for $n\geq 1$, see Theorem \ref{no oriented cycles}.
\end{rema}

Given $n\geq 1$ we will compute the dimensions of  $\tau$-Hochschild (co)homology in degree $n$. They depend on the dimensions of Hochschild (co)homology in degrees strictly smaller than $n$. For $n=1$, we recover the formula we have obtained in \cite{CIBILSLANZILOTTAMARCOSSOLOTAR2024} for $\tau$-Hochschild cohomology in degree one. 

We need the following standard result.

\begin{lemm}\label{alternate}
Let
\begin{equation*}
0 \longrightarrow U_0
\stackrel{\delta_1}{\longrightarrow}
U_1
\stackrel{\delta_2}{\longrightarrow}
U_2
\longrightarrow\cdots
\longrightarrow
U_{n-1}
\stackrel{\delta_n}{\longrightarrow}
U_n
\stackrel{\delta_{n+1}}{\longrightarrow}
0
\end{equation*}
be a finite cochain complex of finite dimensional vector spaces. Let $\HH^i$ be its cohomology at $U_i$. We have
$$\sum_{i=0}^{n}(-1)^i\dim_k\ U_i = \sum_{i=0}^{n}(-1)^i\dim_k\ \HH^i.$$
\end{lemm}

\begin{proof}
We set $\delta_{n+1}=0$. For $ 0\leq i\leq n$ we have
$$ \dim_k\ U_{i}=\dim_k\ \Ker \delta_{i+1} + \dim_k\ \Im \delta_{i+1}.$$
Then 
 \begin{align*} \sum_{i=0}^{n} (-1)^i\dim_k\ U_i= \sum_{i=0}^{n} \dim_k\ (-1)^i\Ker \delta_{i+1} + (-1)^i\dim_k\ \Im \delta_{i+1}
 \end{align*}
and the result follows.\qed
  \end{proof}
  \begin{rema}\sf
We record that for a finite chain complex of finite dimensional vector spaces
\begin{equation*}
0 \longrightarrow V_n
\stackrel{\delta'_n}{\longrightarrow}
V_{n-1}
\stackrel{\delta'_{n-1}}{\longrightarrow}
V_{n-2}
\longrightarrow\cdots
\longrightarrow
V_{1}
\stackrel{\delta'_1}{\longrightarrow}
V_0
\longrightarrow
0
\end{equation*}
with homology $\HH_i$ at $V_i$, the result is 
$$\sum_{j=0}^{n}(-1)^{n-j}\dim_k\ V_j = \sum_{i=0}^{n}(-1)^{n-j}\dim_k\ \HH_j.$$
\end{rema}
  \begin{theo}\label{dim tau HH}
 Let $\Lambda =kQ/I$ be a bound quiver algebra and let $X$ be a $\Lambda$-bimodule.
 For $n\geq 1$ we have
 \begin{itemize} 
\small
\item 
$ \displaystyle \dim_k\ \tauXn= 
\\(-1)^n\left(\sum_{i=0}^{n-1} (-1)^{i+1}\dim_k\  \HH^i(\Lambda, X) +\sum_{\substack{i=0\\x,y\in Q_0}}^{n} (-1)^{i}\dim_k yXx \ \dim_k \Tor_i^\Lambda (k_y, {}_x k) \right)=
\\(-1)^n\left(\sum_{i=0}^{n-1} (-1)^{i+1}\dim_k\  \HH^i(\Lambda, X) + \sum_{\substack{i=0\\x,y\in Q_0}}^{n} (-1)^{i}\dim_k yXx\ \dim_k \Ext^i_\Lambda ({}_xk, {}_yk) \right).$
\item 
$\displaystyle 
\dim_k\ \tauhhX =
\\ (-1)^n\left(\sum_{i=0}^{n-1} (-1)^{i+1}\dim_k\  \HH_i(\Lambda, X) + \sum_{\substack{i=0\\x,y\in Q_0}}^{n}(-1)^{i}\dim_kyXx\  \dim_k\Tor_i^\Lambda (k_x,{}_y k)\right)=\\
\\  (-1)^n\left(\sum_{i=0}^{n-1} (-1)^{i+1}\dim_k\  \HH_i(\Lambda, X) + \sum_{\substack{i=0\\x,y\in Q_0}}^{n} (-1)^{i}\dim_kyXx\  \dim_k\Ext^i_\Lambda ({}_yk,{}_xk)\right).
$
 \end{itemize}
 \normalsize

  \end{theo}
 \begin{proof}

 For short we set, according to Theorem \ref{minimal projective resolution Happel} and using Proposition \ref{pedazos de Tor}:
\begin{equation}\label{An}
\begin{split}
 A_i(X) &=  \Hom_{\Lambda-\Lambda}( \Lambda\otimes_E\Tor^\Lambda_i(E,E)\otimes_E\Lambda,X)\\
 & = \Hom_{E-E}( \Tor^\Lambda_i(E,E),X)\\
 &=\bigoplus_{y,x\in Q_0}\Hom_{E-E}(y\Tor^\Lambda_i(E,E)x, yXx)\\
 &=\bigoplus_{y,x\in Q_0}\Hom_k( \Tor_i^\Lambda (k_y,{}_xk), yXx)
 \end{split}
\end{equation}
and
\begin{equation}\label{Bn}
\begin{split}
 B_i(X) &=  X\otimes_{\Lambda-\Lambda}
 \left(\Lambda\otimes_E\Tor^\Lambda_i(E,E)\otimes_E\Lambda\right)\\& = X\otimes_{E-E} \Tor^\Lambda_i(E,E)
 \\&
 =\bigoplus_{y,x\in Q_0} yXx \otimes x\Tor^\Lambda_i(E,E)y\\&
 = \bigoplus_{y,x\in Q_0} yXx \otimes \Tor^\Lambda_i(k_x,{}_yk).
 \end{split}
\end{equation}

The complexes of cochains (\ref{cochains}) and chains (\ref{chains}) which compute respectively $\HH^*(\Lambda, X)$ and $\HH_*(\Lambda, X)$ are
\begin{equation}\label{cochains with Happel}
0 \longrightarrow A_0(X)
\stackrel{\delta_1}{\longrightarrow}
A_1(X)
\stackrel{\delta_2}{\longrightarrow}
\cdots
\longrightarrow
A_{n-1}(X)
\stackrel{\delta_n}{\longrightarrow}
A_n(X)
\longrightarrow
\cdots
\end{equation}
\begin{equation}\label{chains with Happel}
\cdots \longrightarrow B_n(X) \stackrel{\delta'_n}{\longrightarrow} B_{n-1}(X) \longrightarrow \cdots  \stackrel{\delta'_2}{\longrightarrow} B_1(X) \stackrel{\delta'_1}{\longrightarrow}B_0(X) \longrightarrow 0.
\end{equation}

\normalsize

Since  (\ref{cochains with Happel}) and (\ref{chains with Happel}) are obtained by means of the minimal projective resolution of $\Lambda$, we have by Theorem \ref{coker}, for $n\geq 1$
$$\tauXn= \Coker\ \delta_n
\ \mbox{\  and\  }\ 
\tauhhX=\Ker\ \delta'_n.$$

 Consider the finite cochain complex
 \begin{align*}
0 \rightarrow A_0(X)
\stackrel{\delta_1}{\rightarrow}
A_1(X)&
\stackrel{\delta_2}{\rightarrow}
\cdots
&\rightarrow
A_{n-1}(X)
\stackrel{\delta_n}{\rightarrow}
A_n(X)
\rightarrow
\tauXn
\rightarrow
0.
\end{align*}
It has zero cohomology in $A_n(X)$ and in $\tauXn$, while its cohomology  in $A_i(X)$ for $0\leq i \leq n-1$ is $\HH^i(\Lambda, X)$. Lemma \ref{alternate} gives
$$\sum_{i=0}^{n}(-1)^i\dim_k\  A_i(X) + (-1)^{n+1}\dim_k\  \tauXn = \sum_{i=0}^{n-1}(-1)^i\dim_k\  \HH^i(\Lambda,X).$$

Consider now the finite chain complex
\begin{align*}
0\to \tauhhX \to B_n(X) \stackrel{\delta'_n}{\to} B_{n-1}(X) \to \cdots  \stackrel{\delta'_2}{\to} B_1(X) \stackrel{\delta'_1}{\to}B_0(X) \to 0.
\end{align*}
We get
$$\dim_k\ \tauhhX + \sum_{j=0}^{n}(-1)^{n+1-j}\dim_k B_{j}(X)=
\sum_{j=0}^{n-1}(-1)^{n+1-j}\dim_kH_{j}(\Lambda, X).$$
\qed
 \end{proof}

\begin{coro}\cite{CIBILSLANZILOTTAMARCOSSOLOTAR2024}\label{grado uno}
Let $\Lambda = kQ/I$ be a bound quiver algebra.  We have
$$\dim_k\ \HHH^1_\tau(\Lambda)= \dim_k \HH^0(\Lambda, \Lambda) - \sum_{x\in Q_0}\dim_k x\Lambda x + \sum_{a\in Q_1}\dim_k t(a)\Lambda s(a).$$
$$\dim_k\ \HHH_1^\tau(\Lambda) =  \dim_k \HH_0(\Lambda, \Lambda) - \sum_{x\in Q_0} \dim_k x\Lambda x + \sum_{a\in Q_1} \dim_k s(a)\Lambda t(a).$$
\end{coro}
\begin{proof}
It is well known that  $\Ext^1_\Lambda ({}_xk, {}_yk)$ has a basis in bijection with the set of arrows $a$ such that $s(a)=x$ and $t(a)=y$. On the other hand, $\Ext^0_\Lambda ({}_xk, {}_yk)=\Hom_\Lambda ({}_xk, {}_yk)=0$ if $x\neq y$, and $k$ otherwise.\qed
\end{proof}
To recover precisely the result of \cite{CIBILSLANZILOTTAMARCOSSOLOTAR2024}, note that $Z(\Lambda)= \HH^0 (\Lambda, \Lambda).$ 
\vskip4mm
The formula for local algebras is as follows.

\begin{coro}
Let $\Lambda =kQ/I$ be a local bound quiver algebra, \emph{i.e.} $Q$ has one vertex. Let $X$ be a $\Lambda$-bimodule.
 For $n\geq 1$ we have
 
\begin{align*}
\displaystyle 
 & \dim_k\ \tauXn - \dim_k\ \tauhhX = \\
  &(-1)^n\left(\sum_{i=0}^{n-1} (-1)^{i+1}\dim_k\  H^i(\Lambda, X) - \sum_{i=0}^{n-1} (-1)^{i+1}\dim_k\  H_i(\Lambda, X) \right).
\end{align*}
\end{coro}

Next we give an example showing that in general the $\tau$-Hochschild cohomology and homology are not derived invariant.

\begin{exam}\label{tau not derived}
Let $Q$ be the quiver
\[\begin{tikzcd}
	{{}_x\bullet } & {\bullet_y} & {\bullet_z}
	\arrow["c", shift left, from=1-1, to=1-2]
	\arrow["b", shift left, from=1-2, to=1-1]
	\arrow["d", shift left, from=1-2, to=1-3]
	\arrow["a", shift left, from=1-3, to=1-2]
\end{tikzcd}\]
$I= \langle ada, dc, ad-cb\rangle$ and $\Lambda=kQ/I$.
Let $Q'$ be the quiver 
\[\begin{tikzcd}
	{{}_x\bullet} && {\bullet_y} \\
	& {\stackrel{\bullet}{_z}}
	\arrow["c", from=1-1, to=1-3]
	\arrow["a", from=1-3, to=2-2]
	\arrow["b", from=2-2, to=1-1]
\end{tikzcd}\]
$I'= \langle acba, cbac\rangle$ and $\Lambda'=kQ'/I'$.

The algebras $\Lambda$ and $\Lambda'$ are derived equivalent, see \cite[Example 4.25]{XI 2018}.  
The set $\{1, da, cb+bc\}$ (resp. $\{1, cba, acb+bac\}$) is a basis of the center of $\Lambda$ (resp. $\Lambda'$). Therefore
  $$\dim_k \HHH^0(\Lambda)= 3 =\dim_k \HHH^0(\Lambda').$$
Also,
 $$\dim_k \HHH_0(\Lambda)= 6  = \dim_k \HHH_0(\Lambda').$$
Corollary \ref{grado uno} provides: 
\begin{align*}
&\dim_k\ \HHH^1_\tau(\Lambda) =  3 - (2+2+2) + (1+2+1)= 1,\\
&\dim_k\ \HHH^1_\tau(\Lambda') =  3-(2+2+2)+(1+2+1+1)=2,\\
&\dim_k\  \HHH_1^\tau(\Lambda) =  6 -(2+2+2)+(1+1+1)=3,\\
&\dim_k\  \HHH_1^\tau(\Lambda')  =  6 -(2+2+2)+(1+2+1+1)=5.
\end{align*}
\end{exam}

We will exhibit an example  of a bound quiver algebra of finite global dimension whose $\tau$-homology is non zero:

\begin{exam}\label{no go further}
Let $Q$ be the quiver \[\begin{tikzcd}
	\bullet & \bullet
	\arrow["a", shift left=2, from=1-1, to=1-2]
	\arrow["b", shift left=2, from=1-2, to=1-1]
\end{tikzcd}\] 
and $I=\langle ba\rangle$. The algebra $kQ/I$ is of global dimension $2$. By Corollary \ref{grado uno}, we have 
  $$\dim_k\  \HHH_1^\tau(\Lambda)= 2 - (1+2) + (1+1) = 1. $$
We already know that $\HHH_n^\tau(\Lambda)=0$ for $n\geq 2$ by Proposition \ref{tauHH also} .
\end{exam}

For a bound quiver algebra whose quiver has no oriented cycles, the $\tau$-Hochschild homology vanishes. To prove this result, we need the following well known facts.

\begin{lemm}\label{Bongartz}
Let $\Lambda=kQ/I$ be a bound quiver algebra, and $x,y\in Q_0$. If  there are no paths from $x$ to $y$ in the quiver, namely $y(kQ)x =0$, then for $m\geq 1$ we have
$$\Tor_m^\Lambda (k_y,{}_xk)= 0 =\Ext_\Lambda^m (k_x, k_y).$$  
\end{lemm} 
\begin{proof}
K. Bongartz proves in \cite{BONGARTZ1973} that there are isomorphisms of $E$-bimodules as follows 
$$\Tor_m^\Lambda (E,E)=
\begin{cases}
    \frac{I^n\cap FI^{n-1}F}{FI^n +I^nF} \text{ if } m=2n \text{ for } n\geq 1\\
     \\
     \frac{FI^n\cap I^nF}{I^{n+1} + FI^nF}          \text{ if } m=2n+1 \text{ for } n\geq 0.
\end{cases}
$$
Besides,
$$y\Tor_m^\Lambda (E,E)x= \Tor_m^\Lambda (k_y,{}_xk).$$
Hence 
\begin{align}\label{formula bongartz}
\Tor_m^\Lambda (k_y,{}_xk)=
\begin{cases}
    \frac{y(I^n\cap FI^{n-1}F)x}{y(FI^n +I^nF)x} \text{ if } m=2n \text{ for } n\geq 1\\
     \\
     \frac{y(FI^n\cap I^nF)x}{y(I^{n+1} + FI^nF)x}          \text{ if } m=2n+1 \text{ for } n\geq 0.
\end{cases}
\end{align}
Both numerators are included in $y(kQ)x$ which is $0$. \qed
\end{proof}
\begin{lemm}\label{same vertex}
Let $\Lambda= kQ/I$ be a bound quiver algebra and $x\in Q_0$. If $Q$ has no oriented cycles then $\Tor_m^\Lambda (k_x, {}_xk)=0$ for $m\geq 1$.
\end{lemm}
\begin{proof}
Recall that the ideal $I$ verifies $I\subset F^2$. For $x=y$, the numerators of (\ref{formula bongartz}) are respectively contained in $xF^{2n}x$ for $n\geq 1$ and $xF^{2n+1}x$ for $n\geq 0$. Since $Q$ has no oriented cycles, both are zero.  
 \qed
\end{proof}

\begin{theo}\label{no oriented cycles}
  Let $\Lambda=kQ/I$ be a bound quiver algebra. If $Q$ has no oriented cycles, then 
  $$\HHH_n^\tau(\Lambda)=0 \text{ for all } n\geq 1.$$
\end{theo}
\begin{proof}
Recall that 
 $$B_i(\Lambda) = \bigoplus_{y,x\in Q_0} y\Lambda x \otimes \Tor^\Lambda_i(k_x,{}_yk).$$
 
For $i\geq 1$, we assert that each of the above direct summands vanishes. For the summands with $x=y$, Lemma (\ref{same vertex}) ensures the result.

For $x\neq y$, if $\Tor^\Lambda_i(k_x,{}_yk) = 0$ then the direct summand is zero. If $\Tor^\Lambda_i(k_x,{}_yk)$ is not zero, then $x(kQ)y \neq 0$ by Lemma \ref{Bongartz}. Since $Q$ has no oriented cycles we now that $y(kQ)x=0$, hence $y\Lambda x =0$ and the corresponding direct summand is also zero. 

Hence for $X=\Lambda$ the chain complex (\ref{chains}) is
\begin{align*}
\begin{split}
& \cdots \longrightarrow 0 \stackrel{\delta'_n}{\longrightarrow} 0 \longrightarrow \cdots  \longrightarrow 0 \stackrel{\delta'_1}{\longrightarrow} B_0(\Lambda)=\Lambda\otimes_{\Lambda-\Lambda}P_0 \longrightarrow 0. 
\end{split}
\end{align*}
By Theorem \ref{coker}, $\HHH_n^\tau(\Lambda)=\Ker \delta'_n=0$ for $n\geq 1$.

\qed

\end{proof}

\section{\sf the $\tau$ versions of Happel's question and Han's conjecture}\label{tau Happel's and Han's}

Related to Happel's question and Han's conjecture, we set the following:

\begin{defi}\label{HQ}
Algebras of infinite global dimension with infinite Hochschild cohomology (resp. homology) are called \emph{positive answers to Happel's question (resp.  Han's conjecture).}
\end{defi}

Proposition \ref{tauHH also} leads to the $\tau$ versions of the Happel's question and Han's conjecture.

\begin{defi}\label{tauHQ}
Algebras of infinite global dimension with infinite $\tau$-Hochschild cohomology (resp. homology) are called \emph{positive answers to the $\tau$ version of Happel's question (resp. to  the $\tau$ version of Han's conjecture).}
\end{defi}

Recall from the proof of Theorem \ref{dim tau HH} that  
\begin{itemize}
  \item $ A_n(X) =\Hom_{E-E}(\Tor_n^\Lambda (E,E),X),$
  \item  $B_n(X) = X\otimes_{E-E} \Tor^\Lambda_n(E,E).$
\end{itemize}

 \begin{theo}\label{tau cero implica A cero}
 Let $\Lambda$ be an algebra and $X$ a $\Lambda$-bimodule. Let $N$ be a positive integer.
 We have 
 $$\tauXi=0  \makebox{ for } i\geq N \Leftrightarrow  A_i(X)=0 \  \makebox{ for } i>N.$$
  $$\HH_i^\tau(\Lambda,X)=0 \makebox{ for } i\geq N \Leftrightarrow  B_i(X)=0 \makebox{ for } i>N.$$
 \end{theo}
 \begin{proof}

 Assume $\tauXi=0$ for $i\geq N$. Consider the cochain complex (\ref{cochains with Happel}). By Theorem \ref{coker} we have $\tauXi= \Coker\ \delta_i$, hence our assumption implies that $\delta_i$ is surjective for $i\geq N$.

 Since $\Im \delta_i \subset \Ker\delta_{i+1}$, we infer $\delta_j =0$ for $j>N$, thus $\Coker\ \delta_i = A_i(X)=0$ for $i>N$. 

 The converse is clear. The proof for $\tau$-Hochschild homology is analog.

\qed
 \end{proof}

Next we define two classes of algebras which are clearly of infinite global dimension.

\begin{defi}\label{i+gldim}
 Let $\Lambda = kQ/I$ be a bound quiver algebra. The algebra $\Lambda$ has \emph{infinite +   (resp.  co+) global dimension} if there exists a pair of vertices $(y,x)$ such that
 \begin{itemize}
   \item $y\Lambda x\neq 0$,
   \item $\Tor^\Lambda_*(k_x,{}_yk))$ (resp. $\Tor^\Lambda_*(k_y,{}_xk$) is infinite.
 \end{itemize}
 \end{defi}

For instance non trivial local algebras are of infinite + and infinite  co+ global dimension, see Proposition \ref{local are positive tau H q}.  Note that connected commutative algebras are local, in particular this holds for those algebras. Further examples are shown in the next section.

The following results show that $\tau$-Hochschild cohomology and homology are well adapted to describing the finiteness or not of the previous dimensions.

 \begin{theo}\label{infinite+ is positive for tauHQ}
Let $\Lambda = kQ/I$ be a bound quiver algebra. We have
\begin{itemize}
\item
$\Lambda$ is of infinite co+ global dimension if and only if $\tauL$ is infinite.

\item  $\Lambda$ is of infinite + global dimension if and only if $\HHH_*^\tau(\Lambda)$ is infinite.

\end{itemize}

 \end{theo}
 \begin{proof}
We have
   $A_n(\Lambda) =\oplus_{y,x\in Q_0} \Hom_k(\Tor^\Lambda_n(k_y,{}_xk), y\Lambda x),$
see Proposition \ref{pedazos de Tor}. There exists a pair of vertices $(y,x)$ such that $y\Lambda x\neq 0$ and $\Tor^\Lambda_*(k_y,{}_xk)$ is infinite if and only if $A_*(\Lambda$) is infinite. By Theorem \ref{tau cero implica A cero}, the latter is equivalent  to $\HHH^*_\tau(\Lambda)$ being infinite. The proof of the second statement is analog. \qed

 \end{proof}

\begin{theo}\label{if positive then + and co+}
Let $\Lambda$ be an algebra of infinite global dimension. If $\Lambda$ is a positive answer to Han's conjecture  (resp. to Happel's question), then $\Lambda$ is of infinite + (resp. co+) global dimension.
\end{theo}
\begin{proof}
An algebra verifying the hypothesis is of infinite Hochschild homology (resp. cohomology). The dimensions of $\tau$-Hochschild homology (resp. cohomology) are greater, then  $\tau$-Hochschild homology (resp. cohomology) is infinite as well. The conclusion follows by Theorem \ref{infinite+ is positive for tauHQ}. \qed
\end{proof}
\begin{rema}\label{in other words}
\sf
In other words, we have
\begin{itemize}
\item

An example of an infinite global dimension algebra without being of infinite + global dimension, would be a refutation of Han's conjecture.

\item

An example of an algebra of infinite global dimension but which is nevertheless not of infinite co+ global dimension should be found among the non local algebras which are negative answers to Happel's question. We do not know of any such example.
\item
Assuming that algebras of infinite global dimension are indeed of infinite + global dimension does not directly imply that Han's conjecture is true: if the assumption holds, any algebra of infinite global dimension would have infinite $\tau$-Hochschild homology, by Theorem \ref{infinite+ is positive for tauHQ}. But the dimension of each $\tau$-Hochschild homology space is greater than the dimension of the corresponding Hochschild homology space, meaning that Hochschild homology could still be finite.
\end{itemize}
\end{rema}

\section{\sf Algebras of infinite + and  infinite  co+ global dimension}\label{examples of + and co+}

\subsection{\sf Local algebras}

It is well-known that if a bound quiver algebra $\Lambda=kQ/I$ is local, then $Q$ has a unique vertex. If there are loops, then the algebra is non trivial and is of infinite global dimension.
\begin{prop}\label{local are positive tau H q}
A non trivial local bound quiver algebra $\Lambda$ is of infinite + and  infinite  co+ global dimension. Therefore $\tauhhLs$ and $\tauL$ are infinite.
\end{prop}
\begin{proof}
Let $u$ be the unique vertex of the quiver. We have $u\Lambda u = \Lambda \neq 0$ and $\Tor^\Lambda_*(k_u,{_u}k)$ is infinite. \qed
\end{proof}

\begin{exam}\label{BGMS}
Let $\Lambda_q=k\{x,y\}/\langle x^2, yx+qxy, y^2\rangle$ for $q$ not a root of unity. Note that $\Lambda_q$ is local non trivial, then it is of infinite global dimension. In \cite{BUCHWEITZGREENMADSENSOLBERG2005} it is shown that the algebra $\Lambda_q$ is a negative answer to Happel's question.

More precisely, in degrees $0$, $1$ and $2$ the dimensions of $\HHH^*(\Lambda_q)$ are respectively $2$, $2$ and $1$, see \cite{BUCHWEITZGREENMADSENSOLBERG2005}. In greater degrees it vanishes.

However $\Lambda_q$ is a positive answer to the $\tau$ version of Happel's question, as any non trivial local algebra is, see Proposition \ref{local are positive tau H q}. In the sequel, we compute the  dimensions of the $\tau$-Hochschild cohomology spaces of $\Lambda_q$.

Let $u$ be the unique vertex of the quiver $Q$ with two loops $x$ and $y$. 

Henceforth we will replace $A_n(\Lambda)$ (resp. $B_n(\Lambda)$) by $A_n$ (resp. $B_n$). From \cite{BUCHWEITZGREENMADSENSOLBERG2005} we have
$\dim_k\ \Tor_n^{\Lambda_q}(k_u,{}_uk)= n+1,$ then
$$\dim_k\  A_n = \dim_k\ B_n = (\dim_k\  \Lambda) (n+1) = 4(n+1).$$
Thus
\small
$$\sum_{i=0}^{n}(-1)^i\dim_k\ A_i = \sum_{i=0}^{n}(-1)^i\dim_k\ B_i= 4\sum_{i=0}^{n}(-1)^i (i+1)=
\begin{cases}
               2(n+2) \mbox{ if $n$ is even}\\
               - 2(n+1) \mbox{ if $n$ is odd}\\
            \end{cases}
            $$
            \normalsize
On the other hand, for $n\geq 3$:
$$\sum_{i=0}^{n-1}(-1)^{i+1} \dim_k\ \HHH^i(\Lambda_q) = -2+2-1=-1.$$
According to Theorem \ref{dim tau HH}, for $n\geq 3$ we have
\begin{equation*}
 \dim_k\ \HHH^n_\tau(\Lambda_q) =
            \begin{cases}
             -1 + 2(n+2)& \mbox{ if $n$ is even}\\
              1 + 2(n+1)& \mbox{ if $n$ is odd}
            \end{cases}
 \end{equation*}
Thus for all $n\geq 3$
$$\dim_k\ \HHH^n_\tau(\Lambda_q)=2n+3$$
while
\begin{align*}
  \dim_k\  \HHH^1_\tau(\Lambda_q) &=2-4+8=6\\
  \dim_k\  \HHH^2_\tau(\Lambda_q) &  =-(0-4(1-2+3))=8.
\end{align*}

\end{exam}

\subsection{\sf Finitely generated Yoneda algebras}

Let $\Lambda= kQ/I$ be a bound quiver algebra and $E=\Lambda/r$.
\begin{defi}
  The Yoneda algebra - also called the Ext-algebra - of $\Lambda$ is $\EE(\Lambda)=\Ext_\Lambda^*(E,E)$. Its product is the Yoneda product of exact sequences.
\end{defi}
To each vertex we associate the identity endomorphism of the corresponding simple left module. This way $Q_0$ is a complete system of orthogonal idempotents of $\EE (\Lambda)$.
\begin{prop}\label{idfg}
  Let $\EE$ be a $k$-algebra which is not supposed to be finite dimensional. Assume that $\EE$ is a finitely generated algebra. Let $G_0$ be a complete system of orthogonal idempotents of $\EE$.

   If for every $x\in G_0$ the vector space $x\EE x$ is finite dimensional, then $\EE$ is finite dimensional.
\end{prop}
\begin{proof}
The Peirce decomposition is $\EE = \bigoplus_{ x,y \in G_0} y\EE x.$ The finite system of generators of $\EE$ can be decomposed according to this direct sum, so that we may assume that each element of the finite system of generators lies in a Peirce summand $y\EE x$.

For the purpose of this proof, we define a quiver $G$ associated to this data; its set of vertices vertices is $G_0$. For each generator $g$ in $y\EE x$ there is an arrow $g$ from $x$ to $y$ in $G$. The source (resp. the target) of $g$ is $s(g)=x$ (resp.  $t(g)=y$).

As usual, a path of $G$ is a concatenated sequence of arrows  $\gamma=g_n g_{n-1}\dots g_1$. The sequence of vertices where $\gamma$ \emph{passes through} is $t(g_n), s(g_n), \dots, s(g_1)$.

A path $\gamma=g_n g_{n-1}\dots g_1$ is\emph{ without oriented cycles} if $\gamma$ does not pass through any vertex more than once. Note that any path is a composition of paths without oriented cycles and oriented cycles, which alternate. The number of paths without oriented cycles is finite,  hence the sum below is finite. We have
\begin{align*}
 &\dim_k\ y\EE x \leq &  \\
 &\sum_{\substack{g_n g_{n-1}\dots g_1 \\ \mbox{\footnotesize path from $x$ to $y$}\\\mbox{\footnotesize without } \\\mbox{\footnotesize oriented cycles} }} [\dim_k\  (y\EE y)]\ [\dim_k\  s(g_n)\EE s(g_n)] \dots [\dim_k\  s(g_{2})\EE s(g_{2})]\ [ \dim_k\  x\EE x].
\end{align*}\qed
\end{proof}

Let $\EE(\Lambda)=\Ext^*_\Lambda (E,E)$ be the Yoneda algebra of a bound quiver algebra $\Lambda = kQ/I$. For future use, we consider $\EE(\Lambda)$ as a $k$-category - the Yoneda category of $\Lambda$ - whose set of objects is $Q_0$, morphisms from $x$ to $y$ are $y\EE(\Lambda)x =\Ext^*_\Lambda ({}_xk, {}_yk)$, and composition  is given by the product of the Yoneda algebra.

Recall that the dimension of a $k$-category is the sum of the dimensions of its vector spaces of morphisms.  The Yoneda category $\EE(\Lambda)$ is finite dimensional if and only if $\Lambda$ is of finite global dimension.

\begin{theo}\label{if Yoneda algebra is nice, then + and co+}
Let $\Lambda= kQ/I$ be a bound quiver algebra. Let $\EE(\Lambda)$ be its Yoneda category.
Suppose there exists  a $k$-subcategory $\EE'$ of $\EE(\Lambda)$ which is infinite dimensional although finitely generated. Then $\Lambda$ is of infinite + and  infinite  co+ global dimension, consequently $\tauhhLs$ and $\tauL$ are infinite.
\end{theo}
\begin{proof}
By the previous Proposition \ref{idfg}, there exists $x$ such that
$x\EE ' x$ is infinite dimensional. Therefore $x\EE(\Lambda) x$ is infinite dimensional since $x\EE ' x \subset x\EE(\Lambda) x$. Recall that we have $x\EE x=\Ext^*_\Lambda ({}_xE, {}_xE)$. Of course $x\in x\Lambda x$, hence $x\Lambda x\neq 0$.\qed
\end{proof}

\begin{coro}\label{Koszul}
Let $\Lambda$ be a $n$-Koszul algebra (see for instance \cite{GREEN MARCOS MARTINEZ-VILLA ZHANG}) of infinite global dimension. The algebra $\Lambda$ is of infinite + and  infinite co+ global dimension.
\end{coro}
\begin{proof}
In \cite{GREEN MARCOS MARTINEZ-VILLA ZHANG} it is proven that the Yoneda algebra of $\Lambda$ is generated in degrees $0, 1$ and $2$.\qed
\end{proof}

\begin{exam}
The algebras considered in \cite{PARKER SNASHALL} and \cite{XU ZHANG} are non local negative answers to Happel's question. They are Koszul of infinite global dimension, hence there are of infinite + and  infinite  co+ global dimension by Corollary \ref{Koszul}. Compare with the second item of Remark \ref{in other words}.

\end{exam}

\subsection{\sf Algebras with non zero Peirce components}

The next result generalises the case of local algebras.
\begin{prop}\label{todos no cero}
Let $\Lambda= kQ/I$ be a bound quiver algebra of infinite global dimension. If for each pair of vertices $y,x\in  Q_0$ we have $y\Lambda x \neq 0$ and $x\Lambda y \neq 0$, then $\Lambda$ is of infinite + and  infinite  co+ global dimension.

Also, if for each pair of vertices $y,x\in  Q_0$ we have  $y\Lambda x \neq 0$ and/or $x\Lambda y \neq 0$, then $\Lambda$ is of infinite + and/or  infinite  co+ global dimension.
\end{prop}
\begin{proof}
Consider the decomposition
$$\Tor_*^\Lambda (E,E) = \oplus_{y,x \in Q_0}\Tor_*^\Lambda (k_y, {}_x k).$$
Since $\Tor_*^\Lambda (E,E)$ is infinite, there exist $y,x$ such that $\Tor_*^\Lambda (k_{y}, {}_{x} k)$ is infinite.\qed
\end{proof}

\begin{exam}\label{dos vertices}
Let $Q$ be a quiver with two vertices $x$ and $y$, and let $\Lambda= kQ/I$ be a bound quiver algebra of infinite global dimension. If $Q$  only contains arrows from $x$ to $y$ (or from $y$ to $x$), then $\Lambda$ is of finite global dimension.

Therefore we use Proposition \ref{todos no cero} to infer that $\Lambda$ is of infinite + and  infinite  co+ global dimension.
\end{exam}

\begin{exam}\label{example monomial}
We consider the example  \cite[p. 18]{GREEN ZACHARIA}. Let $Q$ be the quiver 
\[\begin{tikzcd}
	{{}_x\bullet} & {\bullet_y}
	\arrow["a", shift left, from=1-1, to=1-2]
	\arrow["b", shift left=2, from=1-2, to=1-1]
\end{tikzcd}\]
Let $I= \langle aba \rangle$ and $\Lambda=kQ/I$. The graded vector space $\Ext^*_\Lambda ({}_xk, {}_yk)$ is infinite, while all  the other Ext graded vector spaces between simples are finite. Then the Yoneda algebra of $\Lambda$ is not finitely generated. Note that all finitely generated subalgebras of the Yoneda algebra has finite dimension. 
  However $\Lambda$ is of infinite + and  infinite  co+ global dimension, as in Example \ref{dos vertices}.
  \end{exam}

\subsection{\sf Extension conjecture}

For a bound quiver algebra $\Lambda=kQ/I$,  according to \cite{HAPPEL ZACHARIA} the no-loop conjecture was first shown in \cite{LENZING} and reproved in \cite{IGUSA 1990}: if the quiver has a loop, then $\Lambda$ is of infinite global dimension.

The strong no-loop conjecture states that if the quiver has a loop, then the simple module associated to the vertex of the loop is of infinite projective dimension. For $k$ algebraically closed, it has been proved in \cite{IGUSA LIU PAQUETTE}, see also \cite{HAN2015}.

The extension conjecture is as follows, see \cite{HAPPEL ZACHARIA,IGUSA LIU PAQUETTE,LIU MORIN}. If there is a loop at a vertex $u$, then  $\Ext^*({_u}k,{_u}k)$ is infinite - equivalently $\Tor_*^\Lambda(k_u, {_u}k)$ is infinite. Therefore, the following result is clear.

\begin{prop}\label{extension conjecture}
Let $\Lambda$ be a bound quiver algebra such that the quiver contains a loop. If $\Lambda$ verifies the extension conjecture, then $\Lambda$ is of infinite + and  infinite  co+ global dimension. Consequently $\Lambda$ is a positive answer to the $\tau$ version of Han's conjecture and to the $\tau$ version of Happel's question.
\end{prop}

As mentioned in \cite[p. 2741]{IGUSA LIU PAQUETTE}, the extension conjecture is proved for monomial algebras and  special biserial algebras, see \cite{LIU MORIN,GREEN SOLBERG ZACHARIA}. Note that Example \ref{example monomial} is monomial but without loops.

\subsection{\sf Does infinite global dimension imply infinite + or co+ global dimension?}\label{?}

We will make Remark (\ref{in other words}) more precise.

\vskip3mm
Let $\Lambda=kQ/I$ be a bound quiver algebra of infinite global dimension. If $\Lambda$ were not of infinite + global dimension, then
\begin{enumerate}
  \item $\Lambda$ would disprove Han's conjecture, see Theorem \ref{if positive then + and co+},
  \item All the subalgebras of the Yoneda algebra of $\Lambda$ which are infinite dimensional would also be infinitely generated by Theorem \ref{if Yoneda algebra is nice, then + and co+}.
  \item Assume that the extension conjecture is true for the algebra $\Lambda$. Then $Q$ contains no loops by Proposition \ref{extension conjecture}.
\end{enumerate}

Of course we do not know of such an example since, up to date, there are no known  counterexamples to Han's conjecture.
\vskip3mm

Similarly, let $\Lambda$ be a bound quiver algebra of infinite global dimension. If $\Lambda$ were not of infinite co+ global dimension, then:
\begin{itemize}
  \item $\Lambda$ would be a negative answer to Happel's question by Theorem \ref{if positive then + and co+},
      \item Items 2. and 3. would also hold.
\end{itemize}

We do not know of such an example.

\vskip3mm

Finally assume that algebras of infinite global dimension are indeed of infinite + (resp. co+) global dimension. 
Under this assumption, an algebra is of infinite global dimension if and only if its $\tau$-Hochschild homology (resp $\tau$-Hochschild cohomology) is infinite by Theorem \ref{infinite+ is positive for tauHQ}.

\section{\sf Algebras of radical square zero}\label{rad^2=0}

\subsection{\sf Minimal resolution}

For a bound quiver algebra $\Lambda = kQ/I$ there is a well known reduced resolution of $\Lambda$ as $\Lambda$-bimodule.
$$\cdots\rightarrow\Lambda\otimes_E r^{\otimes_E n}\otimes_E\Lambda\stackrel{d_n}{\rightarrow}\cdots\rightarrow\Lambda\otimes_E\Lambda\stackrel{d_0}{\rightarrow}\Lambda\rightarrow 0$$ where the formulas for the differentials are equal to those of the bar resolution.

A bound quiver algebra $\Lambda = kQ/I$ is of \emph{radical square zero} if $I=F^2$, that is all paths of length $2$ are zero in $\Lambda$. In this case $\Lambda=kQ_0 \oplus kQ_1.$ Moreover $r=kQ_1$ and $r^2=0$.

The set of oriented paths of length $n$ of $Q$ is denoted $Q_n$. The vector space with basis $Q_n$ is $kQ_n$.

Actually for radical square zero algebras, the reduced resolution is the minimal one. Indeed, the algebras are monomial and the resolution is Bardzell's one \cite{BARDZELL 1997}. Alternatively, we  clearly have an isomorphism of $E-E$-bimodules  
$r^{\otimes_E n} \equiv kQ_n$ where the sign $\equiv$ means that we consider it as an identification.  For a radical square zero algebra, there are $E-E$-bimodule isomorphisms $\DD\Ext_\Lambda^n(E,E) \simeq kQ_n \simeq \Tor_n^\Lambda(E,E).$  Theorem \ref{minimal projective resolution Happel} ensures that the reduced resolution is the minimal one.

To describe the differentials, we set the following notations.
\begin{itemize}
  \item Let $\gamma = \gamma_n\dots \gamma_1\in Q_n$ where $\gamma_i\in Q_1$ for all $i$. We denote 
  $${}^-\gamma = \gamma_{n-1}\dots \gamma_1 \mbox{\ and \ } \gamma^- = \gamma_n\dots \gamma_2.$$ 
  \item For $a\in Q_1$, we set
  $${}^- a =s(a) \mbox{\ and \ } a^-=t(a).$$
\end{itemize}

\begin{prop}\label{minimal resolution of an algebra radical square zero}
Let $\Lambda=kQ/F^2$ be a radical square zero algebra. The minimal resolution of $\Lambda$ as a $\Lambda$-bimodule is 
\begin{equation}\label{minimal resolution of radical square zero}
\cdots\rightarrow\Lambda\otimes_E kQ_n \otimes_E\Lambda\stackrel{d_n}{\rightarrow}\cdots\rightarrow \Lambda\otimes_E kQ_0\otimes_E\Lambda\stackrel{d_0}{\rightarrow}\Lambda\rightarrow 0
\end{equation}
where $d_n$ for $n\geq 1$ is determined by the morphism of $E-E$-bimodules $$kQ_n\rightarrow \Lambda\otimes_E kQ_{n-1}\otimes_E \Lambda$$ given by
$$\gamma = \gamma_n \dots \gamma_1 \mapsto \gamma_n\otimes {}^- \gamma \otimes s(\gamma) + (-1)^n t(\gamma)\otimes \gamma^- \otimes \gamma_1$$
and $d_0$ is the product of the algebra.
\end{prop}

\subsection{\sf Hochschild and $\tau$-Hochschild homology}
Let $O_n$ be the set of cycles of length $n$, that is $O_n=\{\gamma \in Q_n | s(\gamma)=t(\gamma)\}$. Note that $O_0=Q_0$ and $O_1$ is the set of \emph{loops} of $Q$. The following result is clear.

\begin{lemm}\label{concrete}
Let $\Lambda = kQ/F^2$ and let $X$ be a $\Lambda$-bimodule. We have 
\begin{align*}
  &X\otimes_{\Lambda-\Lambda}\left(\Lambda \otimes_E kQ_n \otimes_E \Lambda\right) & = X\otimes_{E-E} kQ_n  &=\oplus_{y,x \in Q_0} \left(yXx \otimes x(kQ_n) y\right).\\
  &\Lambda\otimes_{\Lambda-\Lambda}\left(\Lambda \otimes_E kQ_n \otimes_E \Lambda\right) & = \Lambda\otimes_{E-E} kQ_n &\equiv kO_n\oplus kO_{n+1}
\end{align*}
where we use $\equiv$ for the following clear identifications: $kQ_0\otimes_{E-E} kQ_n \equiv kO_n$ and $kQ_1\otimes_{E-E} kQ_n \equiv kO_{n+1}$. 
\end{lemm}

\begin{rema}\label{cyclic permutation}
The cyclic group of order $n$ with generator $t$ acts on $O_n$ by cyclic permutations as follows. 
Let $\gamma=\gamma_n\dots \gamma_1$ in $O_n$. Then
  $$t\gamma = \gamma_1\gamma_n\dots \gamma_2 = \gamma_1 \gamma^-.$$ 
We denote $\Omega_n$ the set of orbits of this action. For instance, for $b$ a loop, the number of elements of the orbit of $b^n$ is $1$. 
\end{rema}

\begin{prop}\label{complex of chains for radical square zero}
Let $\Lambda = kQ/F^2$. The chain complex whose homology is  $\HHH_*(\Lambda)$ obtained with the minimal projective resolution (\ref{minimal resolution of radical square zero}) of $\Lambda$ is isomorphic to 
$$\cdots\to kO_n\oplus kO_{n+1}\stackrel{\delta'_n}{\to}kO_{n-1}\oplus kO_n\to \cdots \stackrel{\delta'_1}{\to}kO_0\oplus kO_1\to 0$$
where 
$
\delta'_n = 
\begin{pmatrix}
0 & 0 \\
\mathsf{Id}+ (-1)^nt & 0 
\end{pmatrix}
.$
\end{prop}
\begin{proof}
Let $X$ be a $\Lambda$-bimodule. Using Lemma \ref{concrete}, the boundary map $$X\otimes_{E-E}kQ_n \to X\otimes_{E-E}kQ_{n-1}$$ is the composition
\begin{align*}
&X\otimes_{E-E}kQ_n \longrightarrow \\
&X\otimes_{\Lambda-\Lambda}\left(\Lambda \otimes_E kQ_n \otimes_E \Lambda\right) \stackrel{1_X\otimes d_n}{\longrightarrow} X\otimes_{\Lambda-\Lambda}\left(\Lambda \otimes_E kQ_{n-1} \otimes_E \Lambda\right)\longrightarrow \\
&X\otimes_{E-E}kQ_{n-1}
\end{align*}
which sends
\begin{align*}
&x\otimes \gamma \mapsto\\
&x\otimes t(\gamma)\otimes \gamma_n\otimes \dots \otimes\gamma_1 \otimes s(\gamma) \mapsto\\
& x\otimes t(\gamma)\gamma_n\otimes \gamma_{n-1}\otimes \dots \otimes \gamma_1\otimes s(\gamma) +\\
& (-1)^n x\otimes t(\gamma)\otimes \gamma_n\otimes \gamma_{n-1}\otimes \dots \otimes\gamma_2 \otimes \gamma_1 s(\gamma)=\\
  &x\otimes \gamma_n\otimes \gamma_{n-1}\otimes \dots \otimes \gamma_1\otimes s(\gamma) + \\
  &(-1)^n x\otimes t(\gamma)\otimes \gamma_n\otimes \gamma_{n-1}\otimes \dots \otimes\gamma_2 \otimes \gamma_1 
\mapsto\\ 
&s(\gamma)x\gamma_n\otimes {}^-\gamma +(-1)^n 
\gamma_1xt(\gamma)\otimes \gamma^-=
x\gamma_n\otimes {}^-\gamma +(-1)^n 
\gamma_1x\otimes \gamma^-.
\end{align*}
For $X=\Lambda = kQ_0\oplus kQ_1$ there is an identification 
$$X\otimes_{E-E}kQ_n= (kQ_0\otimes_{E-E}kQ_n) \oplus (kQ_1\otimes_{E-E}kQ_n) \equiv kO_n \oplus kO_{n+1}.$$
Since $r^2=0$, the boundary map restricted to 
$kQ_1\otimes_{E-E}kQ_n\equiv kO_{n+1}$ is zero, while
restricted to $kO_n \equiv (kQ_0\otimes_{E-E}kQ_n)$  its image is contained in $kO_n$  as follows. Let $\gamma \in O_n$ - recall that $s(\gamma)=t(\gamma).$
\begin{align*}
&\gamma \equiv t(\gamma) \otimes \gamma \mapsto \\
& t(\gamma)\gamma_n\otimes {}^-\gamma +(-1)^n 
\gamma_1 t(\gamma) \otimes \gamma^-=\\
&\gamma_n \otimes {}^-\gamma +(-1)^n 
\gamma_1  \otimes \gamma^- \equiv\\
& \gamma +(-1)^n  t\gamma.
\end{align*}\qed
\end{proof}
The next result follows from  Proposition \ref{complex of chains for radical square zero}.
\begin{coro}
Let $\Lambda= kQ/F^2$ be a radical square zero algebra. Its Hochschild homology and cohomology are as follows:  
\begin{align*}
& \HHH_n(\Lambda) \ = \Ker \left(kO_n\ \stackrel{\mathsf{Id} +(-1)^nt}\longrightarrow \ kO_n\right) \ \oplus \ \Coker \left(kO_{n+1}\ \stackrel{\mathsf{Id} +(-1)^{n+1}t}\longrightarrow \ kO_{n+1}\right) \\
&\tauhhL = \Ker \left(kO_n\ \stackrel{\mathsf{Id} +(-1)^nt} \longrightarrow \ kO_n\right)\  \oplus \ kO_{n+1}.
\end{align*}
\end{coro}

We denote $\Omega_n^{\text{even}}$ the set of orbits with an even number of elements. 

\begin{theo}(see \cite[Proposition 3.6]{CIBILS 1990})
Let $\Lambda=kQ/F^2$ be a radical square zero algebra.
\begin{enumerate}[label=(\alph*)]
  \item 
For $n\geq 1$ we have:  
$$\dim_k\ \tauhhL = |\Omega_n| + |O_{n+1}|.$$

\item  If the characteristic of $k$ is different from $2$, then:

$$\dim_k\ \HHH_n(\Lambda)=
\begin{cases}
    |\Omega_n^{\text{even}}| +  |\Omega_{n+1}|& \text{if } n \text{ is even,}\\
    |\Omega_n| +  |\Omega^{\text{even}}_{n+1}|         &\text{if } n \text{ is odd.}
\end{cases}$$

 \item  If the characteristic of $k$ is $2$, then:
 $$\dim_k\ \HHH_n(\Lambda)=
    |\Omega_n| +  |\Omega_{n+1}|.$$
 
\end{enumerate}
\end{theo}   
\begin{proof}
For an orbit $\omega\in\Omega_n$ we denote $k\omega$ the vector space with basis the elements of $\omega$. Then $kO_n= \oplus_{\omega\in\Omega_n}k\omega$, and $\mathsf{Id} +(-1)^nt$ is diagonal with respect of this decomposition.

Assertion 1

$$\dim_k\ \Ker \left(kO_n\ \stackrel{\mathsf{Id} +(-1)^nt}\longrightarrow \ kO_n\right)=|\Omega_n|$$

Let $\omega$ be an orbit of order $b$. Let $\gamma\in\omega$, we have $\omega=\{\gamma, t\gamma, \dots, t^{b-1}\gamma\}$. 
\begin{itemize}
  \item 
  If $n$ is odd, or the characteristic of $k$ is $2$, then for any $n$:  
\begin{align*}
\Ker \left(k\omega\stackrel{\mathsf{Id}  -t}\longrightarrow k\omega\right) =\left\{u\in k\omega\ |\ tu=u\right\}
= k(\gamma+t\gamma+t^2\gamma\dots +t^{b-1}\gamma). 
\end{align*}
  \item 
  If $n$ is even and the characteristic of $k$ is not $2$, then
\begin{align*}\Ker \left(k\omega\stackrel{\mathsf{Id} +t}\longrightarrow k\omega\right) &= \left\{u\in\omega \ |\ tu=-u\right\}
\\&= \begin{cases}
      k(\gamma-t\gamma+t^2\gamma\dots -t^{b-1}\gamma) & \mbox{if } b \text{ is even,} \\
       0 & \mbox{otherwise}.
     \end{cases}    
\end{align*}
\end{itemize}

Assertion 2
 
\begin{itemize}
 \item[-] If the characteristic of $k$ is not $2$, then:
$$\dim_k\ \Coker \left(kO_n\ \stackrel{\mathsf{Id} +(-1)^nt}\longrightarrow \ kO_n\right)=\begin{cases}
    k|\Omega_n|& \text{if } n \text{ is odd,}\\
     k|\Omega_n^{\text{even}}|          & \text{if } n \text{ is even.}
\end{cases}$$
\item[-]   If the characteristic of $k$ is $2$, then:
$$\dim_k\ \Coker \left(kO_n\ \stackrel{\mathsf{Id} +t}\longrightarrow \ kO_n\right)=k|\Omega_n|$$

\end{itemize}

Indeed, let $\omega$ be an orbit of order $b$. Let $\gamma\in\omega$, so that $\omega=\{\gamma, t\gamma, \dots, t^{b-1}\omega\}$.
\begin{itemize}
  \item   
  If $n$ is odd, or in characteristic $2$ for any $n$:
\begin{align*}
&\Coker \left(k\omega\stackrel{\mathsf{Id}  -t}\longrightarrow k\omega\right) = &k\omega/\langle \gamma -t \gamma, \ t\gamma -t^2\gamma, \dots, \ t^{b-1}\gamma - \gamma\rangle = &k\overline{\gamma}.
\end{align*}
  \item 
  If $n$ is even, in characteristic different from $2$:
\begin{align*}
\Coker \left(k\omega\stackrel{\mathsf{Id}  +t}\longrightarrow k\omega\right) &= k\omega/\langle \gamma  + t \gamma,\ t\gamma + t^2\gamma, \dots,\ t^{b-1}\gamma + \gamma\rangle \\&=
\begin{cases}
    k\overline{\gamma}& \text{if } b \text{ is even,}\\
    0           & \text{if } b \text{ is odd.}
\end{cases}
\end{align*}
\end{itemize}
\end{proof}

\subsection{\sf Hochschild and $\tau$-Hochschild cohomology}
This subsection is based on the results of \cite{CIBILS1998}. The computations in \emph{op. cit.} use the minimal resolution of a radical square zero algebra, although this is not mentioned in that paper. Therefore the  computations are also suitable for $\tau$-Hochschild cohomology. We recall some results from \cite{CIBILS1998} which are relevant to us. 

If $U$ and $V$ are sets of paths of a quiver $Q$, we denote 
$$U/\!/V = \{(\gamma, \delta) \in U\times V \ | \ s(\gamma)=s(\delta) \text{ and } t(\gamma)=t(\delta)\}.$$
For instance $Q_n/\!/Q_0 = O_n$, that is the cycles of length $n$. 

\begin{defi}
The linear operator $D_{n+1}:kO_n \longrightarrow k(Q_{n+1}/\!/Q_1)$ is given by 
\begin{equation*}
    D_{n+1}(\gamma) = \sum_{\substack{a\in Q_1\\s(a)=t(\gamma)}} (a\gamma, a)\  + \ (-1)^{n+1} \sum_{\substack{a\in Q_1\\t(a)=s(\gamma)}} (\gamma a, a). 
\end{equation*}
$$D_1(x)=\sum_{\substack{a\in Q_1\\s(a)=x}} (a, a)\  - \  \sum_{\substack{a\in Q_1\\t(a)=x}} (a, a). $$ 
\end{defi}

\begin{theo}\cite[Proposition 2.4]{CIBILS1998}
Let $\Lambda=kQ/F^2$ be a radical square zero algebra. The cochain complex whose cohomology is $\HHH^*(\Lambda)$, given by the minimal resolution (\ref{minimal resolution of radical square zero}) is as follows
\small
$$ 0\to kQ_0\oplus k(Q_0/\!/ Q_1)\stackrel{\delta_1}{\rightarrow} \cdots\to kO_n\oplus k(Q_n/\!/Q_1) \stackrel{\delta_{n+1}}{\rightarrow}kO_{n+1}\oplus k(Q_{n+1}/\!/Q_1)\to \cdots $$
\normalsize
where 
$\delta_{n+1} = 
\begin{pmatrix}
0 & 0 \\
D_{n+1} & 0 
\end{pmatrix}.$

\end{theo}
\begin{coro}For all $n\ge 1$
\begin{align*}
  &\HHH^n(\Lambda) =  \Ker D_{n+1}\oplus\Coker  D_{n}\\
  &\tauLn =  kO_n\oplus\Coker  D_{n} \text{ for } n\geq 1.
\end{align*}
\end{coro}

\begin{defi}\label{crown}
A connected quiver $Q$ is a \emph{$c$-crown} if $Q_0 = \mathbb{Z}/c\mathbb{Z} = Q_1$, where $s: Q_1\to Q_0$ is the identity and $t: Q_1 \to Q_0$ is given by $t(i)=i+1$. 
\end{defi}

\begin{lemm}\cite[Proof of Theorem 3.1]{CIBILS1998}
Let $Q$ be a connected quiver which is not a crown. Let $\Lambda=kQ/F^2$. We have
\begin{itemize}
  \item $D_n$ is injective for $n\geq 2$,
  \item $\Ker D_1 = k\left(\sum_{x\in Q_0}x\right).$
\end{itemize}
\end{lemm}

\begin{theo}\cite[Theorem 3.1]{CIBILS1998}
Let $Q$ be a connected quiver which is not a crown. Let $\Lambda=kQ/F^2$. We have
\begin{itemize}
  \item $\dim_k \HHH^n(\Lambda) = |Q_n/\!/Q_1|-|O_{n-1}|$ for $n\geq 2,$
  \item $\dim_k \HHH^1(\Lambda) = |Q_1/\!/Q_1|-|Q_0|+1,$
  \item $\dim_k \HHH^0(\Lambda) = |Q_0/\!/Q_1|+1.$
\end{itemize}
and 
\begin{itemize}
 \item $\dim_k   \HHH^n_\tau(\Lambda)= |O_n|+|Q_n/\!/Q_1|-|O_{n-1}|$ for $n\geq 2,$
  \item $\dim_k   \HHH^1_\tau(\Lambda)= |O_1|+|Q_1/\!/Q_1|-|Q_0|+1.$
 \end{itemize}
\end{theo}

\footnotesize
\noindent C.C.:\\
Institut Montpelli\'{e}rain Alexander Grothendieck, CNRS, Univ. Montpellier, France.\\
IRL-CNRS  IFUMI-2030\\
{\tt Claude.Cibils@umontpellier.fr}

\medskip
\noindent M.L.:\\
Instituto de Matem\'atica y Estad\'\i stica  ``Rafael Laguardia'', Facultad de Ingenier\'\i a, Universidad de la Rep\'ublica, Uruguay.\\
IRL-CNRS  IFUMI-2030\\
{\tt marclan@fing.edu.uy}

\medskip
\noindent E.N.M.:\\
Departamento de Matem\'atica, IME-USP, Universidade de S\~ao Paulo, Brazil.\\
{\tt enmarcos@ime.usp.br}

\medskip
\noindent A.S.:
\\IMAS-CONICET y Departamento de Matem\'atica,
Facultad de Ciencias Exactas y Naturales,\\
Universidad de Buenos Aires, Argentina. \\{\tt asolotar@dm.uba.ar}

\end{document}